\documentclass[12pt]{article}
\usepackage[marked,printedin%,a4save%
]{myart98}
\ppnum{}{math/9806150}{1998}
\renewcommand{\whereappear}{}
\newcommand{\oper}[1]{\ensuremath{\mathcal{#1}}}
\def\ifundefined#1{\expandafter\ifx\csname#1\endcsname\relax}
\providecommand{\comment}[1]{}

\usepackage{amsfonts}

\newcommand{\algebra}[1]{\ensuremath{\mathfrak{#1}}}
\newcommand{\Heisen}[1]{\ensuremath{\mathbb{H}^{#1}}}

\newcommand{\Cliff}[2][\comment]{\ensuremath{%
\mathcal{C}\kern-0.18em\ell(#1,#2)}}
\newcommand{\object}[2][\,]{\ensuremath{\mathrm{#2}#1}}
\newcommand{\Space}[2]{\ensuremath{\mathbb{#1}^{#2}} }

\newcommand{\FSpace}[2]{\ensuremath{ #1_{#2} }}

\ifundefined{qed}
    \DeclareMathSymbol{\qed}{0}{AMSa}{"03}
\fi
\newcommand{\bra}[1]{\ensuremath{\langle#1\mid}}
\newcommand{\ket}[1]{\ensuremath{\mid#1\rangle}}

\newcommand{\norm}[1]{\left\| #1 \right\|}
\newcommand{\modulus}[1]{\left| #1 \right|}
\newcommand{\scalar}[2]{\left\langle #1,#2 \right\rangle}

\providecommand{\eqref}[1]{\textup{(\ref{#1})}}

\providecommand{\url}[1]{\texttt{#1}}

\def\ip#1#2{#1\cdot#2}

\def\L2{{\FSpace{L}{2}(\Space{R}{n})}}
\def\BB{{\cal B}}

\def\CLie{\Space C{}\mathfrak{h}^n}
\def\Lie{\mathfrak{h}^n}
\newcommand{\SL}{SL(2,\Space{R}{})}

\def\quork#1{{\ifmmode{{\cal C}\kern-0.18em\ell(#1)}\else
${{\cal C}\kern-0.18em\ell(#1)}$\fi}}
\renewcommand\Cliff[1]{\quork{#1}}

\newcommand{\uu}{{\bf u}}
\newcommand{\vv}{{\bf v}}
\renewcommand{\a}{{\bf a}}
\renewcommand{\b}{{\bf b}}
\renewcommand{\c}{{\bf c}}
\newcommand{\p}{{\bf p}}
\newcommand{\q}{{\bf q}}
\newcommand{\w}{{\bf w}}
\newcommand{\x}{{\bf x}}
\newcommand{\y}{{\bf y}}
\newcommand{\z}{{\bf z}}

\newcommand{\SB}{\FSpace{F}{2}( \Space{C}{n} )}

\titleshort{Monogenic Functions and Nilpotent Groups}
\authorshort{J.Cnops and V.Kisil}

\input{cyr.fd}

\begin{document}
\title{Monogenic Functions and\\
Representations of Nilpotent Lie Groups\\
in Quantum Mechanics}
\author{
Jan Cnops\thanks{Post-doctoral researcher, FWO, Belgium.}\\
\normalsize Vakgroep Wiskundige Analyse\\
\normalsize Universiteit Gent\\
\normalsize Galglaan 2, Gent\\
\normalsize B-9000, Belgium\\
\normalsize e-mail: \texttt{jc@cage.rug.ac.be} 
\and
Vladimir V. Kisil\thanks{Suported by the grant 
3GP03196 of the FWO-Vlaanderen (Fund of Scientific Research-Flanders) 
and INTAS grant 93--0322--Ext.}\\
\normalsize Institute of Mathematics\\
\normalsize                Economics and Mechanics\\
\normalsize                Odessa State University\\
\normalsize                ul. Petra Velikogo, 2\\
\normalsize                Odessa, 270057, Ukraine\\
\normalsize                e-mail: \texttt{kisilv@member.ams.org}}
\date{}
\maketitle
\begin{abstract}
We describe several different representations of nilpotent step two Lie 
groups in spaces of monogenic Clifford valued functions. We are inspired by
the classic representation of the Heisenberg group in the Segal-Bargmann
space of holomorphic functions. Connections with quantum mechanics are 
described.
\AMSMSC{46E20}{22E27, 30G35, 30H05, 46E22, 81R05, 81R30}
\end{abstract}
\newpage
\tableofcontents
\section{Introduction}

It is well known, by the celebrated Stone-von Neumann theorem,
that all models for the canonical quantisation~\cite{Mackey63} are
isomorphic and provide us with equivalent representations of the Heisenberg
group~\cite[Chap.~1]{MTaylor86}. Nevertheless it is worthwhile to look
for some models which can act as alternatives for the Schr\"odinger
representation. In particular, the Segal-Bargmann
representation~\cite{Bargmann61,Segal60} serves to 
\begin{itemize}
\item  give a geometric 
representation of the dynamics of harmonic oscillators;
\item present a nice model
for the creation and annihilation operators, which is important for
quantum field theory;
\item allow applying tools of analytic function theory.
\end{itemize}
The huge abilities of the Segal-Bargmann (or Fock~\cite{Fock32}) model are
not yet completely employed, see for example new ideas in a recent
preprint~\cite{NazaikSternin96}.

We look for similar connections between nilpotent Lie groups and spaces of 
monogenic~\cite{BraDelSom82,DelSomSou92} Clifford valued functions. 
Particularly we are interested in a
third possible representation of the Heisenberg group, acting on monogenic
functions on \Space{R}{n}.  There are several reasons why such a model
can be of interest.  First of all the theory of monogenic functions is
(at least) as interesting as several complex variable theory, so the
monogenic model should share many pleasant features with the Segal-Bargmann
model.  Moreover, monogenic functions take their value in a Clifford
algebra, which is a natural environment in which to represent internal
degrees of freedom of elementary particles such as spin.  Thus from the
very beginning it has a structure which in the Segal-Bargmann model has to
be added, usually by means of the second quantization
procedure~\cite{Dirac67}.  So a monogenic representation can be even more
relevant to quantum field theory than the Segal-Bargmann one (see
Remark~\ref{re:quant-field}).

From the different aspects of the Segal-Bargmann space $ \SB $ we select the
one giving a unitary representation of the Heisenberg group
$\Space{H}{n}$.  The
representation is unitary equivalent to the Schr\"odinger representation on
$ \FSpace{L}{2}( \Space{R}{n} ) $ and the Segal-Bargmann transform is
precisely the intertwining operator between these two representations 
(see Appendix~\ref{ss:segal}).
Monogenic functions can be introduced in this
scheme in two ways, as either $ \FSpace{L}{2}( \Space{R}{n} ) $ or $\SB$
can be substituted by a space of monogenic functions. 

In the first case one defines a new unitary irreducible representation of
the Heisenberg group on a space of monogenic functions
and constructs an analogue of the Segal-Bargmann transform as the intertwining
operator of the new representation and the Segal-Bargmann one. We examine this
possibility in section~\ref{se:heisenberg}. In a certain sense the
representation of the Heisenberg group constructed here lies between
Schr\"odinger and Segal-Bargmann ones, taking properties both of them.

In the second case we first select a substitute for the
Heisenberg group, so we can replace the Segal-Bargmann space by a
space of monogenic functions. The space $\Space{C}{n}$
underlying $ \FSpace{F}{2}( \Space{C}{n} ) $ is 
intimately connected with the structure of the Heisenberg group $ 
\Space{H}{n}$ in the sense that
$\Space{C}{n} $ is the quotient of $ \Space{H}{n} $ with respect to its 
centre. In order to define a space of monogenic functions, say on $
\Space{R}{n+1} $, we have to construct a group playing a similar
r\^ole with respect to this space. We describe an option in
section~\ref{se:new}. 

Finally we give the basics of coherent states from square integrable group
representations and an interpretation of the classic Segal--Bargmann space 
in terms of these in Appendix~\ref{se:appendix}.

This paper is closely related to~\cite{Kisil97a}, where connections between
analytic function theories and group representations were described. 
Representations of another group ($\SL$) in spaces of monogenic functions 
can be found in~\cite{Kisil97c}. We hope that the present paper make only
few first steps towards an interesting function theory and other steps will
be done elsewhere.

\section{The Heisenberg group and spaces of analytic functions} 
\label{se:heisenberg}
\subsection{The Schr\"odinger representation of the Heisenberg group}
We recall here some basic facts on the Heisenberg group $\Space{H}{n}$
and its Schr\"odinger representation, see~\cite[Chap.~1]{Folland89}
and~\cite[Chap.~1]{MTaylor86} for details.

The Lie algebra of the Heisenberg group is generated by the $2n+1$ elements
$p_1$, \ldots, $p_n$, $q_1$, \ldots, 
$q_n$, $e$, with the well-known Heisenberg commutator relations:
\begin{equation}\label{eq:a-heisenberg}
[p_i,q_j]=\delta_{ij} e.
\end{equation}
All other commutators vanish. In the standard quantum mechanical
interpretation the operators are momentum and coordinate 
operators~\cite[\S~1.1]{Folland89}.

It is common practice to switch between real and complex Lie algebras.
Complexify $\Lie$ to obtain the complex algebra $\CLie$, and take
four complex numbers  $a$, $b$, $c$ and $d$ such that $ad-bc\not=0$. The
{\it real\/} $2n+1$-dimensional subspace spanned by
$$
A_k=ap_k+bq_k\quad\quad B_k=cp_k+dq_k
$$
and the commutator $[A_k,B_k]=(ad-bc)e$, where $e=[p_k,q_k]$ is of course
isomorphic to $\Lie$, and exponentiating will give a group isomorphic
to the Heisenberg group.

An example of this procedure is obtained from the construction of the
so-called creation and annihilation operators of Bose particles in
the $k$-th state, $a^+_k$ and  $a^{-}_k$ (see~\cite[\S~1.1]{Folland89}).
These are defined by:
\begin{equation}\label{eq:creation}
a^{\pm}_k=\frac{q_k \mp \imath p_k}{\sqrt[]{2}},
\end{equation}
giving the commutators $[a^+_i,a^-_j]=(-\imath)\delta_{ij}e$. Putting
$-\imath e=\ell$, the real algebra spanned by $a^\pm_k$ and $\ell$ is an
alternative realization of $\Lie$, $\Lie_a$.

An element $g$ of the Heisenberg group $\Heisen{n}$ (for any positive
integer $n$) can be represented as $g=(t,\z)$ with $t\in\Space{R}{}$,
$\z=(z_1,\ldots,z_n)\in \Space{C}{n}$. The group law in coordinates
$(t,\z)$ is given by
\begin{equation}\label{eq:g-heisenberg}
g*g'=(t,\z)*(t',\z')=(t+t'+\frac{1}{2}\sum_{j=1}^n\Im(\bar{z}_j
z_j'), \z+\z'),
\end{equation}
where $\Im z$ denotes the imaginary part of the complex number $z$.
Of course the Heisenberg group is non-commutative.

The relation between the Heisenberg group and its Lie algebra is given
by the exponentiation $\exp:\Lie_a\to\Heisen n$. We define the formal
vector $\a^+$ as being $(a^+_1,\ldots,a^+_n)$ and $\a^-$ as
$(a^-_1,\ldots,a^-_n)$, which allows us to use the formal inner products
\begin{eqnarray*}
\uu\cdot \a^+&=&\sum_{k=1}^nu_ka^+_k\\
\vv\cdot \a^-&=&\sum_{k=1}^nv_ka^-_k.
\end{eqnarray*}
With these we define, for real vectors $\uu$ and $\vv$, and real $s$
\begin{eqnarray}
\exp(\uu\cdot(\a^++\a^-))&=&(0,\sqrt 2\uu)\\
\exp(\vv\cdot (\a^--\a^+)&=&(0,\imath\vv)\\
\exp(s \ell)&=&(e^{-2s},0).
\end{eqnarray}

Possible Schr\"odinger representations are parameterized by the
non-zero real number $\hbar$ (the Planck constant). As usual, for
considerations where the correspondence principle between classic and
quantum mechanics is irrelevant, we consider only the case 
$\hbar=1$. The Hilbert space for the Schr\"odinger representation is
$\FSpace{L}{2} (\Space{R}{n})$, where elements of the complex Lie
algebra $\CLie$ are represented by the unbounded operators
\begin{equation} \label{eq:rep-schro}
 \sigma(a^{\pm}_k)=
\frac{1}{\sqrt[]{2}}\left(x_k I \mp \frac{\partial
}{\partial x_k}\right).
\end{equation}
From which it follows, using any $j$, that
$$
\sigma(\ell)=[a^+_j,a^-_j]=-2 I.
$$
The corresponding representation $\pi$ of the Heisenberg group is
given by exponentiation of the $\sigma(a^+_k)$ and $\sigma(a^-_k)$,
but this is most readily expressed by using $p_k$ and $q_k$, and so is
generated by shifts and multiplications $s_\c: f(\x) \mapsto f(\x+\c)$ and
$m_\b: f(\x)\mapsto e^{\imath\ip\x\b}f(\x)$, with the Weyl commutation
relation
\begin{displaymath}
s_\c m_\b = e^{\imath\ip\c\b} m_\b s_\c.
\end{displaymath}

There is an orthonormal basis of $\FSpace{L}{2}(\Space{R}{n})$ on
which the operators $\sigma(a^\pm_k)$ act in an especially simple way. It
consists of the functions:
\begin{equation}\label{eq:hbasis}
\phi_m(\y)=[2^m m! \,\, \sqrt[]{\pi}]^{-1/2} e^{-\ip\x\x/2}H_m(\y),
\end{equation}
where $\y=(y_1,\ldots,y_n)$, $m=(m_1,\ldots,m_n)$, and $H_m(\y)$ is the
generalized Hermite polynomial
$$
H_m(\y)=\prod_{i=1}^nH_{m_i}(y_i).
$$
For these
\begin{displaymath}
a^{+}_k \phi_m(\y)=\sqrt{m_k+1}\,\phi_{m'}(\y), \qquad a^{-}_k \phi_m(\y)=
\sqrt{m_k}\,\phi_{m''}(\y)
\end{displaymath}
where 
\begin{eqnarray*}
m'&=&(m_1,m_2,\ldots, m_{k-1}, m_k + 1, m_{k+1},\ldots,
m_{n})\\
m''&=&(m_1,m_2,\ldots, m_{k-1}, m_k - 1, m_{k+1},\ldots,
m_{n}).
\end{eqnarray*}
This is the most straightforward way to express the creation or
annihilation of a particle in the $k$-th state.

Let us now consider the generating function of the $\phi_m(\x)$,
\begin{equation}\label{eq:sum}
A(\x,\y)=\sum_{j=0}^\infty \frac{x^j}{\sqrt[]{j!}} \phi_k(\y)= \exp(-
\frac{1}{2}(\ip\x\x+\ip\y\y) + \sqrt[]{2}\ip\x\y).
\end{equation}
We state the following elementary fact in Dirac's bra-ket notation.
\begin{lem}\label{le:dirac}
Let $H$ and $H'$ be two Hilbert spaces with orthonormal bases
$\{\phi_k\}$ and $\{\phi'_k\}$ respectively. Then the sum
\begin{equation}\label{eq:operator}
U=\sum_{j=0}^\infty \ket{\phi'_j} \bra{\phi_j}
\end{equation}
defines a unitary operator $U: H \rightarrow H'$ with the following
properties:
\begin{enumerate}
\item $U \phi_k = \phi_k'$;
\item If an operator $T: H\rightarrow H$ is expressed, relative to the
basis $\phi_k$, by the matrix $(a_{ij})$ then the operator $UTU^{-1}: H'
\rightarrow H'$ is expressed relative to the basis $\phi_k'$ by the same
matrix. 
\end{enumerate}
\end{lem}

Now, if we take the function $A(\x,\y)$ from~\eqref{eq:sum} as a kernel
for an integral transform,
\begin{displaymath}
[Af](\y)=\int_{\Space{R}{n}} A(\y,\x) f(\x)\,dx
\end{displaymath}
we can consider it subject to the Lemma above. However, for this we need
to define the space $H'$ and an orthonormal basis $\{\phi_k'\}$  (we
already identified $H$ with $\FSpace{L}{2}(\Space{R}{n})$ and the
$\{\phi_k\}$ are given by~\eqref{eq:hbasis}). There is some freedom in
doing this.

For example it is possible to take the holomorphic extension $A(\z,\y)$
of $A(\x,\y)$ with respect to the first variable.
Then
\begin{enumerate}
\item $H'$ is the Segal-Bargmann space of analytic functions over
$\Space{C}{n} $ with scalar product defined by the integral with respect to
Gaussian measure $e^{-\modulus{\z}^2}\,d\z$;
\item The Heisenberg group acts on the Segal-Bargmann space as follows:
\begin{equation} \label{eq:rep-barg}
[\beta_{(t,\z)}f](\uu)=f(\uu+\z)e^{\imath t-\ip{\bar{\z}}\uu-\modulus{\z}^2/2}.
\end{equation}
This action generates the set of coherent states 
$f_{(0,\vv)}(\uu)=e^{-\bar{\vv}\uu-\modulus{\vv}^2/2}$, $\uu$,
$\vv\in \Space{C}{n}$ from the vacuum vector $f_0(\uu)\equiv 1$;
\item The operators of creation and annihilation are $a^+_k=z_k I$,
$a^-_k=\frac{\partial }{\partial z_k}$.
\item The Segal-Bargmann space is spanned by the
orthonormal basis $\phi_k'=\frac{1}{\sqrt[]{m!}}z^n$ or by the set of
coherent states $f_{(0,\vv)}(\uu)=e^{-\bar{\vv}\uu-\modulus{\vv}^2/2}$, $\uu$,
$\vv\in \Space{C}{n}$
\item The intertwining kernel for $\sigma_{(t,\z)}$~\eqref{eq:rep-schro} and
$\beta_{(t,\z)}$~\eqref{eq:rep-barg} is
\begin{displaymath}
A(\z,\y)=e^{-(\ip\z\z+\ip\x\x)/2-\sqrt[]{2}\ip\z\x}
=\sum_{k=0}^\infty \frac{\z^m}{\sqrt[]{m!}} \cdot
\frac{1}{\sqrt{2^m m!} \sqrt[4]{\pi}} e^{-\ip\x\x/2}H_m(\y)
\end{displaymath}
\item The Segal-Bargmann space has a reproducing kernel
\begin{displaymath}
K(\uu,\vv)= e^{\ip\uu{\bar{\vv}}}=\sum_{k=1}^\infty \phi_k(\uu)
\bar{\phi}_k(\vv)= \int
e^{\ip\uu{\bar{\z}}} e^{\ip\z{\bar{\vv}}} e^{-\modulus{\z}^2}\,d\z.
\end{displaymath}
\end{enumerate}
Details can be found in~\cite{Segal60,Bargmann61}, see also 
Appendix~\ref{ss:segal}.

The Segal-Bargmann space is an interesting and important object, but
there are also other options. In particular we can consider an
alternative representation of the Heisenberg group, this time acting on
monogenic functions, an action we introduce in the next subparagraph.

\subsection{Representation of $\Space{H}{n}$ in spaces of monogenic 
functions}\label{ss:cl-al}
We consider the real Clifford algebra \Cliff{n}, i.e. the algebra
generated by $e_0=1$, ${e_j}$, ${1\leq j \leq n}$, using the identities:
\begin{displaymath}
e_i e_j + e_j e_i = -2 \delta_{ij}, \qquad 1\leq i,j\leq n.
\end{displaymath}
For a function $f$ with values in $\Cliff n$, the action of the Dirac
operator of $\Space{R}{n+1}$ is defined by (here $x=x_0+\x$ is the $n+1$
dimensional variable)
$$
Df(x)=\sum_{i=0}^n\partial_if(x).
$$
A function $f$ satisfying $Df=0$ in a certain domain is called monogenic
there; later on we shall use the term `monogenic' for solutions of more
general Dirac operators. Obviously the notion of monogenicity is closely
related to the one of holomorphy on the complex plane. As a matter of
fact $D^2=-\Delta$, and monogenic functions are solutions of the
Laplacian. The Clifford algebra is not commutative, and so it is
necessary to introduce a symmetrized product. For $k$ elements $a_i$,
$1\leq i\leq k$ of the algebra it is defined by
$$
a_1\times a_2\times\ldots\times a_k=
{1\over k!}\sum_\sigma a_{\sigma(1)}a_{\sigma(2)}\ldots a_{\sigma(n)},
$$
where the sum is taken over all possible permutations of $k$ elements.
If the same element appears several times, we use an exponent notation,
e.g.\ $a^2\times b^3=a\times a\times b\times b\times b$.

Let now $V_k$ be the symmetric power monomial defined by the expression
\begin{equation}
V_k(\x)=\frac{1}{\sqrt[]{k!}}(e_1 x_0-
e_0 x_1)^{k_1} \times (e_2 x_0- e_0 x_2)^{k_2} \times \cdots \times
(e_n x_0- e_0 x_n)^{k_n}.
\end{equation}
It can be proved that these monomials are all monogenic (see e.g.\
\cite{Malonek93}), and even that they constitute a basis for the space of
monogenic polynomials (as a module over $\Cliff n$). In general the
symmetrized product is not associative, and manipulating it can become
quite formal. However, if we restrict the monomials defined above to the
hyperplane $x_0=0$, we obtain
$$
V_k(x)=\frac{1}{\sqrt[]{k!}}x_1^{k_1}x_2^{k_2}\ldots x_n^{k_n},
$$
and so we have the multiplicative property
$$
\sqrt{k!k'!\over(k+k)!}V_kV_{k'}=V_{k+k'},\quad x_0=0.
$$
Another important function is the monogenic exponential function which
is defined by
$$
E(\uu,x)=\exp(\uu\cdot\x)\left(\cos(\norm\uu x_0)-
\frac\uu{\norm\uu}\sin(\uu x_0)\right).
$$
It is not hard to check~\cite[\S~14]{BraDelSom82} that this function is
monogenic, and of course its restriction to the hyperplane $x_0=0$ is
simply the exponential function, $E(\uu,\x)=\exp(\uu\cdot\x)$.

We can therefore extend the symmetric product by the so-called
Cauchy-Kovalevskaya product~\cite[\S~14]{BraDelSom82}: If $f$ and $g$ are
monogenic in $\Space R{n+1}$,
then $f\times g$ is the monogenic function equal to $fg$ on $\Space R n$.
Introducing the monogenic functions $\x_{i}=e_ix_0-e_0x_i$ we can then
write 
$$
V_k(x)=\frac{1}{\sqrt{k!}}x_{1}^{k_1}
\times x_{2}^{k_2}\times\ldots\times x_{n}^{k_n}.
$$

It is fairly easy to check the $V_k$ form an orthonormal set with respect
to the following inner product (see~\cite[\S~3.1]{Cnops94a} on Clifford 
valued inner products):
\begin{equation} \label{eq:inner-m2}
\scalar{V_k}{V_{k'}}=\int_{\Space{R}{n+1}} \bar{V}_k( x) V_{k'}
( x)\, e^{- \modulus{x}^2 } \, dx.
\end{equation}
Let $ \FSpace{M}{2}$ be closure of the linear span of $\{V_k\}$, using
complex coefficients.

The creation and annihilation operators $a^+_k$ and $a^-_k$ can be
represented by symmetric multiplication  (see~\cite{Malonek93}) with the
monogenic variable $\x_j$, which will be written $\x_k I_\times
$, and by the (classical) partial derivative $\frac{\partial }{\partial
\x_j}=\frac{\partial }{\partial {x}_j}$ with respect to $\x_j$,
which appear in the definition of hypercomplex differentiability. On basis
elements they act as follows:
\begin{eqnarray*}
\x_jI_\times
V_{(k_1,\ldots,k_j,\ldots,k_n)}&=&
\sqrt{k_j+1}V_{(k_1,\ldots,k_j+1,\ldots,k_n)},\\
\frac{\partial }{\partial \x_j}
V_{(k_1,\ldots,k_j,\ldots,k_n)}&=&\sqrt{k_j} 
V_{(k_1,\ldots,k_j-1,\ldots,k_n)},
\end{eqnarray*}

It can be checked that this really is a representation of $a^\pm_k$, and
that $a^+_k$ and $a^-_k$ are each other's adjoint. We use the equalities
$a^-_j=\frac{1}{ \sqrt[]{2} }(a^+_j + a^-_j) $ and $a^+_j=
\frac{\imath}{ \sqrt 2}(a^-_j - a^+_j)$, and the commutation relations
$[a^+_i,a^-_j]=e\delta_{ij}$ to obtain a representation of the
Heisenberg group. Thus an element $(t,\z)$, $\z=\uu+i\vv$ of the
Heisenberg group can be written as
\begin{eqnarray*}
(t,\z)&=&
\left(t+\frac{\uu\cdot\uu-\vv\cdot \vv}4,0\right)
\left(0,\frac{(1+\imath)(\uu+\vv)}2\right)
\left(0,\frac{(1-\imath)(\uu-\vv)}2\right)\\
&=&
\exp\left(\left(t+\frac{\uu^2-\vv^2}4\right)e\right)
\exp\left(\frac{(\uu+\vv)q}{\sqrt 2}\right)
\exp\left(\frac{(\uu-\vv)\imath p}{\sqrt 2}\right).
\end{eqnarray*}
It is therefore represented by the operator
\begin{eqnarray}
\pi_{(t,\z)}&=&
\exp\left(-\left(t+\frac{\uu\cdot\uu-\vv\cdot\vv}4\right)\right)\nonumber\\
&&\exp\left(\frac{((\uu+\vv)\cdot\x) I_\times }{\sqrt 2}\right)
\exp\left(\frac{(\uu-\vv)\cdot(\partial_\x)}{\sqrt 2}\right),
\label{eq:rep-m2}
\end{eqnarray}
where obviously for a monogenic function $f$ we have
\begin{eqnarray*}
\exp\left(\frac{(\uu-\vv)\imath p}{\sqrt 2}\right)f(x)&=&
 f\left(x+\frac{\uu-\vv)}{\sqrt 2}\right)\\
\exp\left(\frac{((\uu+\vv)\cdot\x) I_\times }{\sqrt 2}\right)f(x)&=&
E\left(\frac{\uu+\vv}{\sqrt 2},\cdot\right)\times f(x)
\end{eqnarray*}
Therefore it is easy to calculate the image of
the constant function $f_0( \x)= V_0(\x) \equiv 1$, and we obtain
the set of functions
\begin{eqnarray}
f_{(t,\z)}(\x)&=&\pi_{(t,\z)} f_0(\x) \nonumber\\
&=&
\exp\left(-\left(t+\frac{\uu\cdot\uu-\vv\cdot\vv}4\right)\right)
E\left(\frac{\uu+\vv}{\sqrt 2},\cdot\right)\times f_0(x)\nonumber\\
&&\exp\left(-\left(t+\frac{\uu\cdot\uu-\vv\cdot\vv}4\right)\right)
E\left(\frac{\uu+\vv}{\sqrt 2},x\right).
\label{eq:m-coher}
\end{eqnarray}
In the language of quantum physics $f_0(\x)$ is the \emph{vacuum
vector} and functions $f_{(t,\z)}(\x)$ are \emph{coherent states} (or
\emph{wavelets}) for the representation of $ \Space{H}{n}$ we described.
We can summarize the properties of the representation:
\begin{enumerate}
\item  All functions in $\FSpace{M}{2} $ are complex-vector
valued, monogenic in $\Space{R}{n+1}$, and square integrable with
respect to the measure $ e^{ - \modulus{x}^2 }dx$.
\item The representation of the Heisenberg group is given by 
\eqref{eq:rep-m2}. This representation generates a set of coherent
states $f_{(0,z)}(\x)$~\eqref{eq:m-coher} as shifts of the vacuum
vector $f_0(\x) \equiv 1$.
\item The creation and annihilation operators $a^+_k$ and $a^-_k$ are 
represented by  symmetric (Cauchy-Kovalevskaya) multiplication by $\x_j$
and by derivation of monogenic functions. They are adjoint with respect to 
the inner product~\eqref{eq:inner-m2}.
\item $ \FSpace{M}{2} $ is generated as a closed linear space by
the orthonormal basis $V_k(\x)=\frac{1}{\sqrt[]{k!}}(e_1 x_0-
e_0 x_1)^{k_1} \times (e_2 x_0- e_0 x_2)^{k_2} \times \cdots \times
(e_n x_0- e_0 x_n)^{k_n}$, and also by the set of coherent states
$f_{(t,\z)}(\x)$ of~\eqref{eq:m-coher}.
\item The kernel of the operator intertwining the model constructed here
and the Segal-Bargmann one is given by 
\begin{displaymath}
B(\z,x)\sum_{j=0}^\infty V_j(x)
\frac{\z^j}{\sqrt[]{j!}}
=\exp(\sum_{k=1}^n \x_k \bar{z}_k),
\end{displaymath}
which is the holomorphic extension in $\z=\uu+\imath\vv$ of
$E(\uu,x)$. The transformation pair is given by
\begin{eqnarray*}
\BB f(x)&=&\int_{\Space Cn}B(\z,x)f(\z)
\exp\left({-|\z|^2\over 2}\right)\, d\z\\
\BB^{-1}\phi(\z)&=&\int_{\Space R{n+1}}\overline{B(\z,x)}\phi(x)
\exp\left({-|x|^2\over 2}\right)\, dx
\end{eqnarray*}
\item The space $\FSpace{M}{2}$ has a reproducing kernel
\begin{displaymath}
K(x,y)=\sum_{k=0}^\infty V_k(x)
\bar{V}_k(y)=\int_{\Space{C}{n}}B(\z,x)\overline{B(\z,y)}
\, e^{-\modulus{z}^2}dz.
\end{displaymath}
Notice that $\overline{K(x,y)}$ is monogenic in $y$; it is the monogenic
extension of $\overline{E(\y,x)}$.
\end{enumerate}

One can see that some properties of $ \FSpace{M}{2}$ are closer to those
of the Segal-Bargmann space than to those of the space $\FSpace{L}{2}(
\Space{R}{n} )$ it replaces. 
It should be noted that the representation of the Heisenberg
group we obtained here is new and quite unexpected.
\begin{rem} \label{re:quant-field}
We construct $\FSpace{M}{2}$ as a space of complex-vector valued functions.
We can also consider an extended space $\FSpace{\widetilde{M}}{2}$ being
generated by the orthonormal basis $V_k(\x)$ or coherent states
$f_{(0,z)}(\x)$ with Clifford valued coefficients multiplied from the
right hand side. Such a space will share many properties of $\Space{M}{2}$
and have an additional structure: there is a natural representation 
$s: f(\x) \mapsto s^* f(s\x s^*) s$ of $\object{Spin}(n)$ group in
$\FSpace{\widetilde{M}}{2}$. Thus this space provides us with 
a representation of two main symmetries in quantum field theory: the 
Heisenberg group of quantized coordinate and momentum (external degrees of
freedom) and $\object{Spin}(n)$ group of quantified inner degrees of
freedom. Another composition of the Heisenberg group and Clifford algebras 
can be found in~\cite{Kisil93c}.
\end{rem}

\section{Another nilpotent Lie group and its representation} \label{se:new}
\subsection{Clifford algebra and complex vectors}\label{ss:cl-alv}

Starting from the real Clifford algebra \Cliff{n}, we consider complex
$n$-vector valued functions defined on the real line \Space{R}{1} with
values in \Space{C}{n}.  Moreover we will look at the $j$-th component of
\Space{C}{n} as being spanned by the elements $1$ and $e_j$ of the Clifford
algebra.  For two vectors $\uu=(u_1,\ldots,u_n)$ and
$\vv=(v_1,\ldots,v_n)$ we introduce the Clifford vector valued product
(see~\cite[\S~3.1]{Cnops94a} on Clifford valued inner products):
\begin{equation}\label{eq:v-pt}
\uu\cdot \vv=\sum_{j=1}^n \bar{u}_j v_j=\sum_{j=1}^n (u'_j-u''_je_j)
(v'_j+v''_je_j),
\end{equation}
where $u_j=u'_j+u''_j e_j$ and $v_j=v'_j+v''_j e_j$.
Of course, $u\cdot u$ coincides with $\norm{u}^2=\sum_1^n (u_j'^2 +
u_j''^2)$, the standard norm in \Space{C}{n}. So we can introduce the space
$\FSpace{R}{2}(\Space{R}{1})$ of \Space{C}{n}-valued functions on the
real line with the product
\begin{equation}\label{eq:f-pt}
\scalar{f}{f'}=\int_\Space{R}{1} f(x) \cdot f'(x) \, dx.
\end{equation}
Again $\scalar{f}{f}^{1/2}$ gives us the standard norm in the Hilbert space
of $\FSpace{L}{2} $ integrable \Space{C}{n} valued functions.

\subsection{A nilpotent Lie group}
We introduce a nilpotent Lie group, \Space{G}{n}. As a $C^\infty$-manifold
it coincides with \Space{R}{2n+1}. Its Lie algebra has generators $P$,
$Q_j$, $T_j$, $1\leq j\leq n$. The non-trivial commutators between them
are
\begin{equation}
[P,Q_j]=T_j;
\end{equation}
all others vanish.
Particularly \Space{G}{n} is a step two nilpotent Lie group and the
$T_j$ span its centre. It is easy to see that $\Space{G}{1} $ is just
the Heisenberg group $\Space{H}{1} $.

We denote a point $g$ of \Space{G}{n} by $2n+1$-tuple of reals
$(t_1,\ldots,t_n;p;q_1,\ldots,q_n)$. These are the exponential coordinates
corresponding to the basis of the Lie algebra $T_1$, \ldots, $T_n$, $P$,
$Q_1$, \ldots, $Q_n$. The group law is given in exponential coordinates by
the formula
\begin{eqnarray}\label{eq:G-mult}
\lefteqn{ (t_1,\ldots,t_n;p;q_1,\ldots,q_n)*
(t'_1,\ldots,t'_n;p';q'_1,\ldots,q'_n)=}\nonumber\\
&=&(t_1+t_1'+\frac{1}{2}(p'q_1-pq'_1),\ldots,t_n+t_n'+
\frac{1}{2}(p'q_n-pq'_n);\nonumber\\
&&\qquad p+p';q_1+q_1',\ldots,q_n+q_n').
\end{eqnarray}

We consider the homogeneous space $\Omega=\Space{G}{n}/\Space{Z}{}$. Here
$\Space{Z}{}$ is the centre of $\Space{G}{n}$; its Lie algebra is
spanned by $T_j$, $1\leq j \leq n$. It is easy to see that
$\Omega\sim\Space{R}{n+1}$. We define the mapping $s:\Omega\rightarrow
\Space{G}{n} $ by the rule
\begin{equation}
s: (a_0,a_1,\ldots,a_n) \mapsto (0,\ldots,0;
a_0;a_1,\ldots,a_n).
\end{equation}
It is the ``inverse'' of the natural projection $s^{-1}:\Space{G}{n}
\rightarrow \Omega=\Space{G}{n}/\Space{Z}{}$.

It easy to see that the mapping $\Omega\times\Omega \rightarrow \Omega$
defined by the rule $s^{-1}(s(a)*s(a'))$ is just Euclidean
(coordinate-wise) addition $a+a'$.

To introduce the Dirac operator we will need the following set of 
left-invariant differential operators, which generate right shifts on the 
group:
\begin{eqnarray}
T_j &=& \frac{ \partial }{ \partial t_j }, \label{eq:l-first} \\
P&=& \frac{ \partial }{ \partial p } +\frac 12 \sum_1^n q_j \frac{ \partial }{
\partial t_j }, \\
Q_j &=& -\frac{ \partial }{ \partial q_j } + \frac 12 p \frac{ \partial }{
\partial t_j }. \label{eq:l-last}
\end{eqnarray}
The corresponding set of right invariant vector fields generating left 
shifts is
\begin{eqnarray}
T^*_j&=& \frac{ \partial }{ \partial t_j }, \label{eq:r-first} \\
P^*&=& \frac{ \partial }{ \partial p } - \frac 12 \sum_1^n q_j \frac{ \partial }{
\partial t_j }, \\
Q^*_j&=&- \frac{ \partial }{ \partial q_j } - \frac 12 p \frac{ \partial }{
\partial t_j }. \label{eq:r-last}
\end{eqnarray}
A general property is that any left invariant operator commutes with any 
right invariant one.

\subsection{A representation of \Space{G}{n} }
We introduce a representation $\rho$ of \Space{G}{n} in the space
$\FSpace{R}{2}(\Space{R}{})$ by the formula:
\begin{equation}\label{eq:G-rep}
[\rho_g f](x)=(e^{e_1(2t_1+q_1(\sqrt{2}x-p))}f_1(x-\sqrt{2}p), \ldots,
e^{e_n(2t_n+q_n(\sqrt{2}x-p))}f_n(x-\sqrt{2}p)),
\end{equation}
where $f(x)=(f_1(x),\ldots,f_n(x))$ and the meaning of
$\FSpace{R}{2}(\Space{R}{})$ was discussed in Subsection~\ref{ss:cl-alv}.
We note that the generators $e_j$ of Clifford algebras do not interact
with each other under the representation just defined.
One can check directly that~\eqref{eq:G-rep} defines a representation
of \Space{G}{n}. Indeed:
\begin{eqnarray}
%\lefteqn
[\rho_g\rho_{g'}f](x)
&=&\rho_g(e^{e_1(2t'_1+q'_1(\sqrt{2}x-p'))}f_1(x-\sqrt{2}p'), \ldots,
\nonumber \\
&&\qquad \qquad e^{e_n(2t'_n+q'_n(\sqrt{2}x-p'))}f_n(x-\sqrt{2}p')) 
\nonumber \\
&=& (e^{e_1(2t_1+q_1(\sqrt{2}x-p))}
e^{e_1(2t'_1+q'_1(\sqrt{2}(x-\sqrt{2}p)-p'))}
f_1(x-\sqrt{2}p-\sqrt{2}p'), \nonumber \\
&& \ldots, \nonumber \\
&& e^{e_n(2t_n+q_n(\sqrt{2}x-p))}
e^{e_n(2t'_n+q'_n(\sqrt{2}(x-\sqrt{2}p)-p'))}
f_n(x-\sqrt{2}p-\sqrt{2}p')) \nonumber \\
&=& (e^{e_1(2(t_1+t'_1+\frac{1}{2}(p'q_1-pq'_1))+
(q_1+q'_1)(\sqrt{2}x-(p+p')))}
f_1(x-\sqrt{2}(p+p')), \nonumber\\
&& \ldots, \nonumber \\
&& (e^{e_n(2(t_n+t'_n+\frac{1}{2}(p'q_n-pq'_n))+
(q_n+q'_n)(\sqrt{2}x-(p+p')))}
f_n(x-\sqrt{2}(p+p')) \nonumber \\
&=& [\rho_{gg'}f](x),
\end{eqnarray}
where $gg'$ is defined by~\eqref{eq:G-mult}.

$\rho_g$ has the important property that it preserves
the product~\eqref{eq:f-pt}. Indeed:
\begin{eqnarray}
\scalar{\rho_gf}{\rho_gf'}&=& \int_\Space{R}{} [\rho_g f](x) \cdot
[\rho_gf'](x)  \,dx \nonumber \\
&=& \int_\Space{R}{}\sum_{j=1}^n
 \bar{f}_j(x-\sqrt[]{2}p)e^{-e_j(2t_j+q_j(\sqrt[]{2}x-p))} \nonumber \\
&& \qquad\quad\quad e^{e_j(2t_j+q_j(\sqrt[]{2}x-p))} f'_j(x-\sqrt[]{2}p)
 \,dx \nonumber \\
&=& \int_\Space{R}{}\sum_{j=1}^n
\bar{f}_j(x-\sqrt[]{2}p)
f'_j(x-\sqrt[]{2}p)  \,dx \nonumber \\
&=& \int_\Space{R}{}\sum_{j=1}^n
 \bar{f}_j(x)  f'_j(x)  \,dx \nonumber \\
&=&\scalar{f}{f'}. \nonumber
\end{eqnarray}
Thus $\rho_g$ is \emph{unitary} with respect to the Clifford valued inner
product~\eqref{eq:f-pt}. Notice this notion is stronger than unitarity for the scalar
valued inner product, as the latter is the trace of the Clifford valued
one. A proof of unitarity could also consist of proving the action of the
Lie algebra is skew-symmetric, i.e.\ that for an element $b$ of the Lie
algebra and $f$ arbitrary
$$
\scalar{d\rho_bf}{f}=\scalar{f}{-d\rho_bf}.
$$
Here $d\rho_b$ is derived representation of $d$ for an element $b\in 
\algebra{g}_n$ of the Lie algebra of $\Space{G}{n}$. In the next subsection
we will need the explicit form of it. For the selected basis of 
$\algebra{g}_n$ we have:
\begin{eqnarray} 
[d\rho(T_j) f] (x) & = &  (0, 0,\ldots, 0,2e_1 f_j(x),0, \ldots ,0,0);
\nonumber \\
{} [d\rho (P) f] (x) & = & (- \sqrt[]{2} \frac{\partial}{\partial x}
f_1(x), \ldots, -\sqrt[]{2} \frac{\partial}{\partial x} f_j(x), \ldots , -
\sqrt[]{2} \frac{\partial}{\partial x}
f_n(x)); \nonumber \\
{} [d\rho (Q_j) f] (x) & = & ( 0,0, \ldots,0, \sqrt[]{2} e_j x f_j(x), 0,
\ldots , 0,0).
\nonumber
\end{eqnarray}
Particularly $d\rho(Q_j) d\rho( Q_k)=0$ for all $j \neq k$. This does not 
follow from the structure of $\Space{G}{}$ but is a feature of the 
described representation.
\begin{rem}
The group $\Space{G}{n}$ is called as ``a generalized Heisenberg group'' 
in~\cite{Kumahara97} where its induced representations are considered.
\end{rem}

\subsection{The wavelet transform for \Space{G}{n} }
In $\FSpace{R}{2}(\Space{R}{})$ we have the \Space{C}{n}-valued function
\begin{equation}
f_0(x)=(e^{-x^2/2},\ldots,e^{-x^2/2}),
\end{equation}
which which is the \emph{vacuum vector} in this case.  It is a zero
eigenvector for the operator
\begin{equation}
\label{eq:g-annihil}
a^-=d\rho(P)- \sum_{j=1}^n e_j d\rho(Q_j),
\end{equation}
which is \emph{the only annihilating operator} in this model. But we still
have $n$ creation operators:
\begin{equation}
a_k^+ = d\rho(P)- \sum_{j=1}^n (1-2\delta_{jk})e_j d\rho(Q_j)= a^{-} + 
2e_kd\rho(Q_k).
\end{equation}
While $a^-$ and $a_k^+$ look a little bit exotic for 
$\Space{G}{1}=\Space{H}{1}$ they are exactly the standard annihilation and 
creation operators.
Another feature of the representation is that the $a_k^+$ do not commute
with each other and have a non-trivial commutator with $a^-$:
\begin{displaymath}
[a_j^+, a_k^+]= 2e_kd\rho(T_k) -2e_j d\rho(T_j), \qquad 
[a_j^+,a^-]=-2e_jd\rho(T_j)
\end{displaymath}

We need the transforms of $f_0(x)$ under the action~\eqref{eq:G-rep}, i.e.
the \emph{coherent states} $f_g(x)=[\rho_gf_0](x)$ in this model:
\begin{eqnarray}
f_g(x)&=&(\ldots,
e^{e_j(2t_j+q_j(\sqrt[]{2}x-p))} e^{-(x-\sqrt[]{2}p)^2/2},\ldots)
\nonumber \\
&=&(\ldots,
e^{2e_jr_j-(p^2+q_j^2)/2} e^{-((p-e_j q_j)^2+x^2)/2+
\sqrt[]{2}(p-e_jq_j)x},\ldots) \nonumber \\
&=&(\ldots,
e^{2e_jt_j-z_j\bar{z}_j/2} e^{-(\bar{z}_j^2+x^2)/2+
\sqrt[]{2}\bar{z}_jx},\ldots) \nonumber \\
\end{eqnarray}
where $z_j=p+e_j q_j$, $\bar{z}_j=p-e_j q_j$.

Having defined coherent states we can introduce the \emph{wavelet
transform} $\oper{W}: \FSpace{R}{2}(\Space{R}{}) \rightarrow 
\FSpace{L}{\infty}(\Space{G}{n})$ by the standard formula:
\begin{equation} \label{eq:g-wt}
 \oper{W} f (g) = \scalar f{f_g}.
\end{equation}

Calculations completely analogous to those of the complex case allow us to
find the images $\oper Wf_{(t',\a)}(t,\z)$ of coherent states
$f_{(t',\a)}(x)$ under~\eqref{eq:g-wt} as follows:
\begin{eqnarray}
\oper Wf_{(t',\a)}(t,\z)&=&\scalar{f_{(t',\a)}}{f_{(t,\z)}} \nonumber\\
&=&\int_\Space{R}{} \sum_{j=1}^n
\exp\left(-2e_jt_j-\frac{z_j\bar{z}_j}2-\frac{{z}_j^2+x^2}2+
\sqrt[]{2}z_jx\right)\nonumber\\
&&\quad
\exp\left(+2e_jt'_j-\frac{a_j\bar{a}_j}2-\frac{\bar{a}_j^2+x^2}2
+\sqrt[]{2}\bar{a}_jx\right)\,dx \nonumber \\
&=&
\sum_{j=1}^n
\exp\left(-2e_jt_j-\frac{z_j\bar{z}_j}2+2e_jt'_j-\frac{a_j\bar{a}_j}2
+\bar{a}_j z_j\right)\nonumber\\
&&\qquad\times \int_\Space{R}{} \exp({-x^2+
2x\frac{z_j+\bar{a}_j}{\sqrt[]{2}}-\frac{(z_j+\bar{a}_j)^2}{2}})\,dx
\nonumber \\
&=&
\sum_{j=1}^n
\exp\left(-2e_j(t_j-t'_j)-\frac{z_j\bar{z}_j+a_j\bar{a}_j}2
+\bar{a}_j z_j\right)\nonumber\\
&&\qquad\times \int_\Space{R}{}
\exp(-(x-\frac{z_j+\bar{a}_j}{\sqrt{2}})^2)\,dx \nonumber \\
&=&
\sum_{j=1}^n
\exp\left(-2e_j(t_j-t'_j)-\frac{z_j\bar{z}_j+a_j\bar{a}_j}2
+\bar{a}_j z_j\right)
\end{eqnarray}
Here $a_j=a_0+e_j a_j$, $\bar{a}_j=a_0-e_ja_j$; $z_j$, $\bar{z}_j$ were
defined above. 

In this case all $\oper Wf_{(t',\a)}(t,\z)$ are \emph{monogenic}
functions with respect to the following Dirac operator:
\begin{equation} \label{eq:dirac-f}
\frac{\partial}{\partial p}-\sum_{j=1}^n e_j \frac{\partial}{\partial q_j}
+\frac 12\sum_{j=1}^n  (e_j p+ q_j)\frac{\partial}{\partial t_j},
\end{equation}
with $z_j$ related to $p$ and $q_j$ as above. This can be checked by the
direct calculation or follows from the observation: the 
Dirac operator~\eqref{eq:dirac-f} is the image of the annihilation operator 
$a^-$~\eqref{eq:g-annihil} under the wavelet transform~\eqref{eq:g-wt}. 
The situation is completely analogous to the Segal-Bargmann case, where
holomorphy is defined by the operators $\frac{\partial}{\partial
\bar{z}_k}$, which are the images of the annihilation operators $a^-_k$.
Actually, it is the Dirac operator associated with the unique left invariant
metric on $\Space Gn/\Space Z{}$ for which $P$ together with the $Q_k$
forms an orthonormal basis in the origin, and therefore everywhere.

The operator~\eqref{eq:dirac-f} is a realization of a generic Dirac
operator constructed for a nilpotent Lie group, see~\cite{ConMosc82}.
Indeed the operator~\eqref{eq:dirac-f} is defined by the formula $D= P + 
\sum_1^n e_j Q_j$, where $P$ and $Q_j$  are the left invariant vector fields 
in \eqref{eq:l-first}--\eqref{eq:l-last}. So the operator~\eqref{eq:dirac-f}
is left invariant and one has only to check the monogenicity of $ 
\oper Wf_{(0,0)}(t,\z) $---all other functions
$\oper Wf_{(t',\a)}(t,\z) $ are its left shifts.

Of course all linear combinations of the $\oper Wf_{(t',\a)}(t,\z)$ are
also monogenic.
So if we define two function spaces, $\FSpace{R}{2}$ and
$\FSpace{M}{2}$, as being the closure of the linear span of all $f_g(x)$ and
$\oper Wf_{(t',\a)}(t,\z)$ respectively, then
\begin{enumerate}
\item $\FSpace{M}{2}$ is a space of monogenic function on \Space{G}{n} in
the sense above.
\item \Space{G}{n} has representations both in $\FSpace{R}{2}$ and in
$\FSpace{M}{2}$. On the second space the group acts via left regular
representation.
\item These representation are intertwining by the integral
transformation with the kernel $T(t',\a,x)=f_{(t',\a)}(x)$.
\item The space $\FSpace{M}{2}$ has a reproducing kernel
$K(t',\a,t,z)=\oper Wf_{(t',\a)}(t,\z)$.
\end{enumerate}
The standard wavelet transform can be processed as expected. 

For the reduced wavelet transform associated with the mapping 
$s:\Omega \rightarrow \Space{G}{n}$ in particular we have
$$
\widehat{\oper W}f_\a(\z)=\oper Wf_{(0,\a)}(0,\z)=\sum_{j=1}^n\exp\bar{a}_j z_j.
$$

However the reduced wavelet transform cannot be constructed from a
single vacuum vector.  We need exactly $n$ linearly independent vacuum
vectors and the corresponding multiresolution wavelet analysis (wavelet
transform with several independent vacuum vectors) which is outlined
in~\cite{BratJorg97a} (see also M.G.~Krein's works~\cite{Krein48a} on
``directing functionals").  Indeed we have $n$ different vacuum vectors
$(\ldots, 0, e^{-x^2/2},0,\dots)$ each of which is an eigenfunction for the
action of the centre of $ \Space{G}{n}$.  All functions $\widehat{\oper
W}f_\a(\z)$ are \emph{monogenic} with respect to the Dirac operator
\begin{equation}
D=\frac{\partial}{\partial p}+ \sum_{j=1}^n e_j \frac{\partial}{\partial q_j}.
\end{equation}
For details on the reduced wavelet transform in a more general setting we
refer to the second Appendix.

\setcounter{section}{0}
\renewcommand{\thesection}{\Alph{section}}
\section{Appendices} \label{se:appendix}
\subsection{The wavelet transform and coherent states}
Let $X$ be a topological space and $G$ be a group of transformations
$g: x \mapsto g \cdot x$ acting from
the left on $X$, i.e.\ $g_1 \cdot(g_2 \cdot x)=(g_1 g_2)\cdot 
x$. Moreover, assume $G$ acts transitively on $X$. Let there
exist a measure $dx$ on $ X$ and  a representation $\pi_g: f(x)
\mapsto m(g,x) f( g^{-1}\cdot x)$ (where $m(g,x)$ is a function), such that
$\pi$ is unitary with respect to the scalar product
$\scalar{f_1}{f_2}_{\FSpace{L}{2}(X ) } =
\int_X f_1(x) \bar{f}_2(x) \,dx$, i.e.
\begin{displaymath} 
\scalar{\pi_gf_1}{\pi_gf_2}_{\FSpace{L}{2}( X ) }
= \scalar{f_1}{f_2}_{\FSpace{L}{2}( X ) }\qquad \forall f_1,
f_2 \in \FSpace{L}{2}(X) .
\end{displaymath}
We shall work with the Hilbert space $\FSpace{L}{2}( X )$ where each
$\pi_g$ is a unitary action.

Let $H$ be a closed compact\footnote{While the compactness will be 
explicitly used during our abstract consideration, it is not crucial 
in fact. Appendix~\ref{ss:segal} will show how to deal with
non-compact $H$.} subgroup of $G$ and let $f_0(x)$ be a function on which
each element $h$ of $H$ acts as multiplication with a constant $\chi(h)$,
\begin{equation} \label{eq:homogenious}
\pi_hf_0(x)=\chi(h) f_0(x) \qquad , \forall h\in H.
\end{equation}
This means $\chi$ is a character of $H$ and $f_0$ is a common
eigenfunction for all operators $\pi_h$. Equivalently $f_0$ is a
common eigenfunction for the operators 
corresponding under $\pi$ to a basis of the Lie algebra of $H$.
Note also that $ \modulus{\chi(h)}^2=1 $ because $\pi$ is unitary.
$f_0$ is called \emph{vacuum vector} (with respect to the subgroup $H$).
We introduce $\FSpace{F}{2}(X)$, the closed linear subspace of 
$\FSpace{L}{2}(X)$ uniquely defined by the conditions:
\begin{enumerate}
\item\label{it:begin} $f_0\in \FSpace{F}{2}( X )$;
\item $\FSpace{F}{2}( X )$ is $G$-invariant;
\item\label{it:end} $\FSpace{F}{2}( X )$ is 
$G$-irreducible.
\end{enumerate}
$f_0$ is then called a {\it cyclic vector} for this space. \ref{it:end}\
puts an extra condition upon $f_0$: there could be functions $f_0$ the
orbit of which spans a reducible space. The theory can be extended to
this case without much difficulty, but we will restrict ourselves to
irreducible spaces here, and the restriction of $\pi$ on $
\FSpace{F}{2}(X) $ is an irreducible unitary representation. The
transforms of $f_0$ will be called \emph{coherent states}. They will be
written down as $w_g$, with
$$
w_g(x)=f_0(g^{-1}\cdot x).
$$

The \emph{wavelet transform} $\oper{W}$ can be defined 
for square-integrable unitary representations $\pi$ by
the formula~\cite{Kisil95a} 
\begin{eqnarray}
\oper{W}&:& \FSpace{F}{2}( X ) \rightarrow
\FSpace{L}{ \infty }(G) \nonumber \\
&:& f(x) \mapsto
\oper Wf(g)=\scalar f{w_g}_{\FSpace{L}{2}( X ) } \label{eq:wavelets}
\end{eqnarray}
The main advantage of the wavelet transform $\oper{W}$ is that it 
expresses $\pi$ in geometrical terms in the sense that it
\emph{intertwines} $\pi$ and the left regular representation $\lambda$ on
$G$ defined by $\lambda_gF(g')=F(g^{-1}g')$:
\begin{equation} \label{eq:g-inter}
\lambda_g\oper{W} f(g')
=\oper{W}f(g^{-1}g')
= \scalar{f}{w_{g^{-1}g'}}
= \scalar{\pi_g f}{w_{g'}}
=\oper{W}\pi_g f(g'),
\end{equation}
i.e., $\lambda \oper{W} = \oper{W} \pi$. Applying this for $f=f_0$ gives
\begin{equation} \label{eq:g-winter}
\oper{W}f_0(g^{-1}g')
=\oper{W}w_g(g').
\end{equation}
Another important feature of $\oper W$ is that it does not lose
information: the function $ f $ can \emph{be recovered} as a linear
combination of the coherent states $w_g$ from its wavelet transform
$\oper Wf(g)$~\cite{Kisil95a}:
\begin{equation} \label{eq:g-inverse}
f(x)=\int_G \oper Wf(g) w_g(x)\, dg
\,dg.
\end{equation}
Here $dg$ is the Haar measure on $G$ which is normalized in such a way that
$\int_G \modulus{ \oper Wf_0(g)}^2\,dg=1 $. One also has the
orthogonal \emph{projection} $\oper P$ from $
\FSpace{L}{2}(G,dg)$ onto the image $ \FSpace{F}{2}(G,dg)=\oper W
\FSpace{F}{2}(X) $, which is just the
convolution on $G$ with the image $ \oper Wf_0(g)$
of the vacuum vector~\cite{Kisil95a}:
\begin{equation} \label{eq:g-proj}
\oper{P}\phi(g')=\int_{G} \phi(g)  \oper Wf_0(g^{-1}g')\,dg.
\end{equation}

\subsection{The reduced wavelet transform}
Our main observation will be that one can be much more economical (if
the subgroup $H$ is non-trivial) with the help of~\eqref{eq:homogenious}: in
this case one does not need to know $ \oper Wf(g) $ on the whole
group $G$, but only on the homogeneous space $G/H$.

Let $s: G \rightarrow G$ be a mapping such that $[s(b)]=[b]$ (square
brackets denoting equivalence classes in $G/H$), and such
that $s(a)=s(b)$ if $[a]=[b]$. Let $\Omega$ be the image of $G$ under
$s$. Any $g\in G$ has a unique decomposition of the form $g=s(g)h$,
$a\in \Omega$, and we will write $h=r(g)=s(g)^{-1}g$. $G/H$ is a left
$G$-homogeneous space for the action defined by $g: [a] \mapsto [ga]$.
Therefore $\Omega$ can be considered to be a $G$-homogeneous space by
the action $t_g:a\mapsto s(ga)$. Due to~\eqref{eq:homogenious} we have
\begin{eqnarray}\label{eq:character}
w_g(x)&=&\pi_gf_0(x)=\pi_{s(g)}(\pi_{r(g)}f_0)(x)\nonumber\\
&=&\pi_{s(g)}(\chi(r(g))f_0)(x)
=\chi(r(g))\pi_{s(g)}f_0(x)\nonumber\\
&=&\chi(r(g))w_{s(g)}(x).
\end{eqnarray}
Therefore
$$
\oper Wf(g) =\scalar f{w_g}_{\FSpace{L}{2}( X ) }
=\overline\chi(r(g))\scalar f{w_{s(g)}}_{\FSpace{L}{2}( X ) }
=\overline\chi(r(g))\oper Wf(r(g)).
$$
Thus $\oper Wf(g)$ is known once its restriction to $\Omega$ is known or,
in more abstract sense, once the wavelet transform is known on $G/H$.
Therefore the restriction of $\oper Wf$ to $\Omega$ merits a new
notation, $\widehat{\oper W}f$. The mapping $\widehat{\oper W}:
 \FSpace{F}{2}(X) \rightarrow
\FSpace{L}{ \infty }( \Omega ) $ will be called \emph{reduced
wavelet transform} and we shall denote by $ \FSpace{F}{2}(\Omega)$ the
image of $\widehat{\oper W}$ equipped with the inner product induced
by $ \widehat{\oper W}$ from $ \FSpace{F}{2}(X) $.

It follows from~\eqref{eq:g-inter} that $ \widehat{\oper W} $
intertwines $\pi$ with the representation $\rho$ given for a function
$\phi$ on $\Omega$ by
$$
\rho_g\phi(a)=\chi(r(g^{-1}a))\phi(s(g^{-1}a)).
$$
Indeed, for $\phi$ of the form $\phi=\widehat{\oper W}f$ we have that
$$
\rho_g\widehat{\oper W}f(a)=\chi(r(g^{-1}a))\widehat{\oper W}f(s(g^{-1}a))=
\oper Wf(g^{-1}a),
$$
and so
\begin{equation} \label{eq:a-inter}
\rho_g \widehat{\oper W}f(a) = \oper Wf(g^{-1}a)=\widehat{\oper W}\pi_gf(a).
\end{equation}
While $\rho$ is not as geometrical as $\lambda$, in applications 
it is still has a more geometrical nature than the original $\pi$.
If the Haar measure $dh$ on $H$ is taken in such a way that
$\int_H  \modulus{ \chi(h) }^2  \, dh=1$ and $dg=dh\,da$
we can rewrite~\eqref{eq:g-inverse}
as follows:
\begin{eqnarray*}
f(x) & = &  \int_G \oper Wf(g) w_g(x)\, dg \nonumber \\
& = & \int_\Omega \int_H \oper W f(ah)w_{ah}(x)\, dh \, da \nonumber \\
& = & \int_\Omega \int_H \widehat{\oper W} f(a)\overline\chi(h)\chi(h)w_a(x)\,
 dh \, da \nonumber \\
& = & \int_{ \Omega } \widehat{\oper W}f(a) w_a(x)\,da
\end{eqnarray*}
We define
an integral transformation $ \oper{F} $ according to the last formula:
\begin{equation} 
\oper F\phi(x)=\int_\Omega\phi(a)w_a(x)\, da \label{eq:a-inverse}.
\end{equation}
This has the property $ \oper{F}\circ \widehat{\oper W} = I$ on $
\FSpace{F}{2}(X) $. One can then consider the integral transform
$\oper K=\oper F\circ\widehat{\oper W}$,
explicitly
\begin{equation} \label{eq:szego}
\oper Kf(x)=
\int_\Omega \scalar f{w_a}_{\FSpace{L}{2}( X ) }  w_a(x) \, da,
\end{equation}
which is defined on the whole  of $ \FSpace{L}{2}(X)$ (not only $
\FSpace{F}{2}(X) $). It is known that $ \oper K$ \emph{is an
orthogonal projection $ \FSpace{L}{2}(X) \rightarrow \FSpace{F}{2}(X)
$}~\cite{Kisil95a}. If we formally use linearity of the scalar product $
\scalar{\cdot}{\cdot}_{ \FSpace{L}{2}(X)}$ (i.e., assume that
Fubini's Theorem holds) we obtain from~\eqref{eq:szego}

\begin{eqnarray}
\oper Kf(x)
&=&
\int_\Omega \scalar f{w_a}_{\FSpace{L}{2}( X ) }  w_a(x) \, da \nonumber \\
\comment{
& = &  \scalar{f(\y)}{\int_{ \Omega }f_{s(a)}(\y) \bar{f}_{s(a)}(x) \, da
}_{\FSpace{L}{2}(
X )}   \nonumber \\
}
& = & \int_X f(\y) K(y,x)\, d\mu(\y) \label{eq:bergman},
\end{eqnarray}
where
\begin{displaymath}
K(y,x)=\int_{ \Omega } \bar w_a(\y)  w_a(x) \, da
\end{displaymath}
Sometimes a reduced form $ \widehat{\oper{P}}: \FSpace{L}{2}( \Omega )
\rightarrow \FSpace{F}{2}( \Omega )$ of the projection $\oper
P$~\eqref{eq:g-proj} is of interest in itself. It is an extension of the
integral operator $\widehat{\oper W}\circ\oper F$, and it is an easy calculation
using \eqref{eq:g-winter}~that
\begin{eqnarray} \label{eq:s-b-proj}
[\widehat{\oper{P}} \phi](a')=\int_\Omega \phi(a) \oper Wf_0(a^{-1}a')
\bar{\chi}(r(a^{-1}a')) \, da.
\end{eqnarray}
As we shall see its explicit form can be calculated easily in practical
cases.

Observe that, from~\eqref{eq:g-inter}, the image of $\oper{W}$ is
invariant under action of the left but not right regular
representations. $ \FSpace{F}{2}(\Omega)$ is invariant under the
representation~\eqref{eq:a-inter}, which is a pullback of the left
regular representation on $G$, but not its right counterpart, and so in
general there is no way to define an action of left-invariant vector
fields on $\Omega$, which are infinitesimal generators of right
translations, on $ \FSpace{L}{2}(\Omega) $. But there is an exception.
Let $ \algebra{X}_j $ be a maximal set of left-invariant vector fields
on $G$ such that
\begin{displaymath}
\algebra{X}_j \oper Wf_0(g)=0.
\end{displaymath}
Because the $\algebra{X}_j $ are left invariant and \eqref{eq:g-winter}\
we have $\algebra{X}_j \oper Ww_{g'}(g)=0$ for all $g'$ and thus the image
of $\oper{W}$, being the linear span of $\oper Ww_{g'}$, is part of the
intersection of kernels of $\algebra{X}_j$. The same remains true if we
consider the pullback $ \widehat{\algebra{X}}_j$ of $\algebra{X}_j$ to $
\Omega$. Note that in general there are fewer linearly independent
$\widehat{\algebra{X}}_j$ than there are ${\algebra{X}}_j$. We call $
\widehat{\algebra{X}}_j$ \emph{Cauchy-Riemann-Dirac} operators because
of the property that
\begin{equation} \label{eq:dirac}
\widehat{\algebra{X}}_j \widehat{\oper W}f(a)=0 \qquad \forall \oper Wf
\in \FSpace{F}{2}(\Omega).
\end{equation}
Explicit constructions of the Dirac type operator for a discrete 
series representation can be found in 
\cite{AtiyahSchmid80,KnappWallach76}.

\subsection{The Segal-Bargmann space} \label{ss:segal}
We consider a representation of the Heisenberg group $ \Space{H}{n}
$ (see Section~\ref{se:heisenberg}) on $ \FSpace{L}{2}( \Space{R}{n} ) $ by
shift and multiplication operators~\cite[\S~1.1]{MTaylor86}:
\begin{equation} \label{eq:schrodinger}
g=(t,\z): f(\x) \rightarrow [\pi_{(t,\z)}f](\x)=
e^{\imath(2t-\sqrt{2}\ip\q\x+\ip\q\p)}
f(\x- \sqrt 2\p), \qquad \z=\p+\imath q,
\end{equation}
This is the Schr\"odinger representation with parameter $\hbar=1$. As a
subgroup $H$ we select the centre of \Space{H}{n} consisting of elements
$(t,0)$. It is non-compact but using the special form of
representation~\eqref{eq:schrodinger} we can consider the
cosets\footnote{$ \widetilde{G}$ is sometimes called the \emph{reduced}
Heisenberg group. It seems that $ \widetilde{G}$ is a virtual object,
which is important in connection with a selected representation of $G$.}
$ \widetilde{G} $ and $ \widetilde{H} $ of $G$ and $H$ by the subgroup
with elements $(\pi m,0)$, $m\in \Space{Z}{}$.
Then~\eqref{eq:schrodinger} also defines a representation of $
\widetilde{G} $ and $ \widetilde{H} \sim \Gamma $. We consider the Haar
measure on $ \widetilde{G} $ such that its restriction on $
\widetilde{H} $ has total mass equal to $1$.

As ``vacuum vector'' we will select the original \emph{vacuum
vector} of quantum mechanics---the Gauss function $f_0(\x)=e^{-\ip\x\x/2}$.
Its transformations are defined as follows:
\begin{eqnarray*}
w_g(\x)=\pi_{(t,\z)} f_0(\x) & = &
e^{\imath(2t-\sqrt{2}\ip\q\x+\ip\q\p)}
\,e^{-{(\x- \sqrt 2\p)}^2/2}\\
& = & e^{2\imath t-(\ip\p\p+\ip\q\q)/2}
e^{- ((\p-\imath \q)^2+\ip\x\x)/2+\sqrt{2}\ip{(\p-\imath \q)}\x}
\\
& = & e^{2\imath t-\ip\z{\bar\z}/2}e^{- (\ip{\bar\z}{\bar\z}+\ip\x\x)/2
+\sqrt{2}\ip{\bar\z}\x}.
\end{eqnarray*}
In particular $w_{(t,0)}(\x)=e^{-2it}f_0(\x)$, i.e.\ it really is
a vacuum vector with respect to $\widetilde{H}$ in the sense of our
definition. Of course $\widetilde{G} / \widetilde{H}$ is
isomorphic to  $\Space{C}{n}$. Embedding $\Space Cn$ in $G$ by the
identification of $(0,\z)$ with $\z$, the mapping $s:
\widetilde G \rightarrow \widetilde{G}$ is defined simply by
$s((t,\z))=(0,\z)=\z$; $\Omega$ then is identical with $\Space Cn$.

The Haar measure on $ \Space{H}{n} $ coincides with the standard
Lebesgue measure on $ \Space{R}{2n+1} $~\cite[\S~1.1]{MTaylor86} and so
the invariant measure on $ \Omega $ also coincides with Lebesgue
measure on $\Space{C}{n}$. Note also that the composition law sending
$\z_1$ $\z_2$ to $s((0,\z_1)(0,\z_2))$ reduces to Euclidean shifts on $
\Space{C}{n} $. We also find $s((0,\z_1)^{-1}\cdot (0,\z_2))=\z_2-\z_1$
and $r((0,\z_1)^{-1}\cdot (0,\z_2))= (\frac{1}{2} \Im\ip{\bar\z_1}{\z_2},0)$.

The reduced wavelet transform takes the form of a mapping
$\FSpace{L}{2}(\Space{R}{n} ) \rightarrow \FSpace{L}{2}( \Space{C}{n} ) 
$ and is given by the formula
\begin{eqnarray}
\widehat{\oper W}f(\z)&=&\scalar{f}{w_{(0,\z)}}\nonumber \\
     &=&\pi^{-n/4}\int_\Space{R}{n} f(\x)\, e^{-\ip\z{\bar\z}/2}\,e^{-
(\ip\z\z+\ip\x\x)/2+\sqrt{2}\ip\z\x}\,dx \nonumber \\
     &=&e^{-\modulus{\z}^2/2}\pi^{-n/4}\int_\Space{R}{n} f(\x)\,e^{-
(\ip\z\z+\ip\x\x)/2+\sqrt{2}\ip\z\x}\,dx, \label{eq:tr-bargmann}
\end{eqnarray}
where $\z=\p+\imath\q$. Then $\widehat{\oper W}f$ belongs to
$\FSpace{L}{2}( \Space{C}{n} , dg)$. This can better be expressed by
saying that the 
function $\breve{f}(\z)=e^{\modulus{\z}^2/2}\widehat{\oper W}f(\z)$ belongs
to $\FSpace{L}{2}( \Space{C}{n} , e^{- \modulus{\z}^2 }dg)$
because $\breve{f}(\z)$ is analytic in $\z$. These functions constitute the
\emph{Segal-Bargmann space}~\cite{Bargmann61,Segal60}  
$ \FSpace{F}{2}( \Space{C}{n}, e^{-
\modulus{\z}^2 }dg) $ of functions analytic in $\z$ and
square-integrable with respect the Gaussian measure $e^{-
\modulus{\z}^2}d\z$. Analyticity of $\breve{f}(\z)$ is equivalent to
the condition that $( \frac{ \partial }{ \partial\bar{\z}_j } + \frac{1}{2} \z_j
I ) \oper Wf(\z)$ equals zero.

The integral in~\eqref{eq:tr-bargmann} is the well-known
Segal-Bargmann transform~\cite{Bargmann61,Segal60}. Its inverse is 
given by a realization of~\eqref{eq:a-inverse}:
\begin{eqnarray}
f(\x) & = & \int_{ \Space{C}{n} } \widehat{\oper Wf(\z)} w_{(0,\z)}(\x)\,d\z
\nonumber\\
& = & \int_{ \Space{C}{n} } \breve{f}(\z)  e^{-
(\bar{\z}^2+\ip\x\x)/2+\sqrt{2}\bar{\z}x}\, e^{- \modulus{\z}^2}\, d\z.
\end{eqnarray}
This gives~\eqref{eq:a-inverse} the name of Segal-Bargmann 
inverse. The corresponding operator $\oper{P}$~\eqref{eq:szego} is the
identity operator $ \FSpace{L}{2}(\Space{R}{n}) \rightarrow 
\FSpace{L}{2}(\Space{R}{n}) $ and~\eqref{eq:szego} gives an integral
presentation of the Dirac delta. 

Meanwhile the orthoprojection  
$ \FSpace{L}{2}( \Space{C}{n},  e^{- \modulus{\z}^2 }dg)  \rightarrow
\FSpace{F}{2}( \Space{C}{n},  e^{- \modulus{\z}^2 }dg) $ is of interest
and is a principal ingredient in Berezin
quantisation~\cite{Berezin88,Coburn94a}. We can easy find its kernel
from~\eqref{eq:s-b-proj}. Indeed, $ \widehat{\oper
W}f_0(\z)=e^{-\modulus{\z}^2}$, and the kernel is
\begin{eqnarray*}
K(\z,\w) & = & \widehat{\oper W}f_0(\z^{-1}\cdot \w)
\bar{\chi}(r(\z^{-1}\cdot \w))\\
& = & \widehat{\oper W}f_0(\w-\z)\exp(\imath\Im(\ip{\bar\z}\w) \\
& = & \exp(\frac{1}{2}(- \modulus{\w-\z}^2 +\ip\w{\bar\z}-\ip\z{\bar\w}))\\
& = &\exp(\frac{1}{2}(- \modulus{\z}^2- \modulus{\w}^2) +\ip\w{\bar\z}).
\end{eqnarray*}
To obtain the reproducing kernel for functions
$\breve{f}(\z)=e^{\modulus{\z}^2} \widehat{\oper W}f(\z) $ in the
Segal-Bargmann space we multiply $K(\z,\w)$ by $e^{(-\modulus{\z}^2+
\modulus{\w}^2)/2}$ which gives the standard reproducing kernel, $\exp(-
\modulus{\z}^2 +\ip\w{\bar\z})$ \cite[(1.10)]{Bargmann61}.

\renewcommand{\thesection}{}
\section{Acknowledgments}
The paper was written while the second author stayed at the Department
of Mathematical Analysis, University of Gent  whose hospitality and
support he gratefully acknowledges. The stay was provided for by the
grant 3GP03196 of the FWO-Vlaanderen (Fund of Scientific
Research-Flanders), Scientific Research Network ``Fundamental Methods
and Technique in Mathematics'' and INTAS grant 93--0322--Ext
subsequently.

\small
\bibliographystyle{plain}
\bibliography{ABBREVMR,akisil,analyse,aphysics,arare}

\newcommand{\noopsort}[1]{} \newcommand{\printfirst}[2]{#1}
  \newcommand{\singleletter}[1]{#1} \newcommand{\switchargs}[2]{#2#1}
  \newcommand{\irm}{\textup{I}} \newcommand{\iirm}{\textup{II}}
  \newcommand{\vrm}{\textup{V}}
  \providecommand{\htmladdnormallink}[2]{#1}\providecommand{\MR}[1]{\textbf{MR%
}\# #1}
\begin{thebibliography}{10}

\bibitem{AtiyahSchmid80}
Michael Atiyah and Wilfried Schmid.
\newblock A geometric construction of the discrete series for semisimple {Lie}
  group.
\newblock In J.A. Wolf, M.~Cahen, and M.~De Wilde, editors, {\em Harmonic
  Analysis and Representations of Semisimple {Lie} Group}, volume~5 of {\em
  Mathematical Physics and Applied Mathematics}, pages 317--383. D. Reidel
  Publishing Company, Dordrecht, Holland, {\noopsort{}}1980.

\bibitem{Bargmann61}
V.~Bargmann.
\newblock On a {H}ilbert space of analytic functions and an associated integral
  transform. {Part I}.
\newblock {\em Comm. Pure Appl. Math.}, 3:215--228, 1961.

\bibitem{Berezin88}
Felix~A. Berezin.
\newblock {\em Method of Second Quantization}.
\newblock ``Nauka'', Moscow, {\noopsort{}}1988.

\bibitem{Deinze93}
F.~Brackx, R.~Delanghe, and H.~Serras, editors.
\newblock {\em Clifford Algebras and Their Applications in Mathematical
  Physics}, volume~55 of {\em Fundamental Theories of Physics}, Dordrecht,
  1993. Kluwer Academic Publishers Group.
\newblock \MR{94j:00019}.

\bibitem{BraDelSom82}
F.~Brackx, R.~Delanghe, and F.~Sommen.
\newblock {\em Clifford Analysis}, volume~76 of {\em Research Notes in
  Mathematics}.
\newblock Pitman Advanced Publishing Program, Boston, 1982.

\bibitem{BratJorg97a}
Ola Bratteli and Palle E.~T. Jorgensen.
\newblock Isometries, shifts, {C}untz algebras and multiresolution wavelet
  analysis of scale ${N}$.
\newblock {\em Integral Equations Operator Theory}, 28(4):382--443, 1997.
\newblock E-print \texttt{funct-an/9612003}.

\bibitem{Cnops94a}
Jan Cnops.
\newblock {\em {Hurwitz} Pairs and Applications of {M\"obius} Transformations}.
\newblock {Habilitation} dissertation, Universiteit Gent, Faculteit van de
  Wetenschappen, 1994.
\newblock \texttt{ftp://cage.rug.ac.be/pub/clifford/jc9401.tex}.

\bibitem{Coburn94a}
Lewis~A. Coburn.
\newblock {Berezin-Toeplitz} quantization.
\newblock In {\em Algebraic Mettods in Operator Theory}, pages 101--108.
  Birkh\"auser Verlag, New York, 1994.

\bibitem{ConMosc82}
Alain Connes and Henri Moscovici.
\newblock The $l\sp{2}$-index theorem for homogeneous spaces of {Lie} groups.
\newblock {\em Ann. of Math.}, 115(2):291--330, 1982.

\bibitem{DelSomSou92}
Richard Delanghe, Frank Sommen, and Vladimir Sou\v{c}ek.
\newblock {\em Clifford Algebra and Spinor-Valued Functions}.
\newblock Kluwer Academic Publishers, Dordrecht, {\noopsort{}}1992.

\bibitem{Dirac67}
P.~A.~M. Dirac.
\newblock {\em Lectures on Quantum Field Theory}.
\newblock Yeshiva University, New York, {\noopsort{}}1967.

\bibitem{Fock32}
V.~A. Fock.
\newblock Konfigurationsraum und zweite quantelung.
\newblock {\em Z. Phys. A}, 75:622--647, 1932.

\bibitem{Folland89}
Gerald~B. Folland.
\newblock {\em Harmonic Analysis in Phase Space}.
\newblock Princeton University Press, Princeton, New Jersey, {\noopsort{}}1989.

\bibitem{Kisil93c}
Vladimir~V. Kisil.
\newblock Clifford valued convolution operator algebras on the {Heisenberg}
  group. {A} quantum field theory model.
\newblock In Brackx et~al. \cite{Deinze93}, pages 287--294.
\newblock \MR{1266878}.

\bibitem{Kisil95a}
Vladimir~V. Kisil.
\newblock Integral representation and coherent states.
\newblock {\em Bull. Belg. Math. Soc. Simon Stevin}, 2(5):529--540, 1995.
\newblock \MR{97b:22012}.

\bibitem{Kisil97a}
Vladimir~V. Kisil.
\newblock Two approaches to non-commutative geometry.
\newblock page~35, 1997.
\newblock \eprint{funct-an/9703001}{http://xxx.lanl.gov/abs/funct-an/9703001/}.

\bibitem{Kisil97c}
Vladimir~V. Kisil.
\newblock Analysis in{ $\Space{R}{1,1}$} or the principal function theory.
\newblock {\em Complex Variables Theory Appl.}, page~25, 1998.
\newblock (To
  appear)\eprint{funct-an/9712003}{http://xxx.lanl.gov/abs/funct-an/9712003/}.

\bibitem{KnappWallach76}
A.W. Knapp and N.R. Wallach.
\newblock Szeg\"o kernels associated with discrete series.
\newblock {\em Invent. Math.}, 34(2):163--200, 1976.

\bibitem{Krein48a}
M.~G. Kre{\u\i}n.
\newblock On {H}ermitian operators with directed functionals.
\newblock {\em Akad. Nauk Ukrain. RSR. Zbirnik Prac\cprime\ Inst. Mat.},
  1948(10):83--106, 1948.
\newblock \MR{14:56c}, reprinted in~\cite{KreinII}.

\bibitem{KreinII}
M.~G. Kre{\u\i}n.
\newblock {\em {\cyr {I}zbrannye Trudy}. {II}}.
\newblock Akad. Nauk Ukrainy Inst. Mat., Kiev, 1997.
\newblock \MR{96m:01030}.

\bibitem{Kumahara97}
Keisaku Kumahara.
\newblock On non-unitary induced representations of a generalized {Heisenberg}
  group.
\newblock {\em Math. Jap.}, 46(1):5--14, 1997.

\bibitem{Mackey63}
George~W. Mackey.
\newblock {\em Mathematical Foundations of Quantum Mechanics}.
\newblock W.~A.~Benjamin, Inc., New York, Amsterdam, {\noopsort{}}1963.

\bibitem{Malonek93}
Helmuth~R. Malonek.
\newblock Hypercomplex differentiability and its applications.
\newblock In Brackx et~al. \cite{Deinze93}, pages 141--150.
\newblock \MR{94j:00019}.

\bibitem{NazaikSternin96}
Vladimir Nazaikinskii and Boris Sternin.
\newblock Wave packet transform in symplectic geometry and asymptotic
  quantization.
\newblock In Komrakov B.P., Krasil'shchik I.S., Litvinov G.L., and Sossinsky
  A.B., editors, {\em Lie Groups and Lie Algebras. Their Representations,
  Generalizations and Applications}, number 433 in Mathematics and Its
  Applications, pages 47--70, Dordrecht-Boston-London, 1998. Kluwer Academic
  Publishers.

\bibitem{Segal60}
Irving~E. Segal.
\newblock {\em Mathematical Problems of Relativistic Physics}, volume~II of
  {\em Proceedings of the Summer Seminar (Boulder, Colorado, 1960)}.
\newblock American Mathematical Society, Providence, R.I., 1963.

\bibitem{MTaylor86}
Michael~E. Taylor.
\newblock {\em Noncommutative Harmonic Analysis}, volume~22 of {\em Math. Surv.
  and Monographs}.
\newblock American Mathematical Society, Providence, R.I., {\noopsort{}}1986.

\end{thebibliography}
\end{document}
%
%\begin{eqnarray}
%\pi((0,(1+\imath)\uu))=\exp(\sigma(\uu\cdot a^+))&=&(0,\uu)\\
%\exp(\vv\cdot a^-)&=&(0,\imath\vv)\\
%\exp(s e)&=&(e^{-2s},0).
%\end{eqnarray}